# LARGE DEVIATIONS FOR OCCUPATION TIMES OF MARKOV PROCESSES WITH $L_2$ SEMIGROUPS

By Naresh Jain and Nicolai Krylov[1]

*University of Minnesota*

Our aim is to unify and extend the large deviation upper and lower bounds for the occupation times of a Markov process with $L_2$ semigroups under minimal conditions on the state space and the process trajectories; for example, no strong Markov property is needed. The methods used here apply in both continuous and discrete time. We present the proofs for continuous time only because of the inherent technical difficulties in that situation; the proofs can be adapted for discrete time in a straightforward manner.

**1. Introduction.** Let $(E, \mathcal{B}(E), m)$ be a $\sigma$-finite measure space with $m(E) > 0$. $B(E)$ denotes the space of extended real-valued measurable functions on $E$. For $p \geq 1$, let

$$\hat{L}_p(m) = \left\{ f \in B(E) : \int |f|^p \, dm < \infty \right\}.$$

Members of $\hat{L}_p(m)$ are functions. $L_p(m)$ will denote the corresponding space whose members are equivalence classes, as usual. The $L_p$-norm will be denoted by $\| \cdot \|_p$ and for $p = 2$, the inner product in $L_2(m)$ will be denoted by $\langle \cdot, \cdot \rangle$. Our goal is to establish lower and local upper large deviation probability bounds for occupation times of Markov processes with semigroups $\{T_t, t \in T\}$ acting on $L_2(m)$ under minimal conditions, where $T = \{0, 1, 2, \ldots\}$ or $T = [0, \infty)$. In the latter case, we assume that the semigroup is strongly continuous on $L_2(m)$, that is

$$(1.1) \qquad \lim_{t \to 0} \|T_t f - f\|_2 = 0, \qquad f \in L_2(m).$$

There is a large class of Markov processes where the semigroups are strongly continuous on $L_2(m)$. Such processes include most of the Lévy

---

Received September 2006; revised October 2007.

[1]Supported in part by NSF Grants DMS-01-40405 and DMS-06-53121.

processes with $E = \mathbb{R}^n$ and $m$ being Lebesgue measure, as well as many diffusion processes and stationary Markov processes which are ergodic with an invariant probability measure $m$. However, the study of large deviations for occupation times of such processes has been done piecemeal under different sets of conditions in a number of papers. Most notable among these is the fundamental work of Donsker and Varadhan [3] and the paper of Fukushima and Takeda [5]. In [3] and their subsequent work, Donsker and Varadhan define the rate function in terms of the infinitesimal generator acting on bounded measurable functions $B_b(E)$ and in [3], under certain conditions, they identify the rate function in terms of the associated Dirichlet form when the process is symmetric and $m$ is the $\sigma$-finite invariant measure. In [5], Fukushima and Takeda assume $E$ to be a locally compact complete separable metric space and assume the process to be symmetric and ergodic, with $m$ the invariant probability measure. They show that in their setting, the $m$-exceptional set that appears in the upper and lower bounds has capacity 0, where the notion of capacity is the same as for a Hunt process; they also identify the rate function in terms of the associated Dirichlet form. Motivated by [5], Mück [9] generalizes the result to the case where $E$ is a "topological Lusin space," not necessarily metrizable. The existence of the underlying "right process" with a notion of "capacity" is provided by Ma and Röckner [7], under the assumption that the associated Dirichlet form (not necessarily symmetric) is "quasiregular." Related to the work in [3], there is also the interesting work of Deuschel and Stroock [2].

In the current paper, we introduce the notion of an $m$-thin set, which turns out to have capacity 0 in the setting of [9]. Our exceptional sets will be $m$-thin, so the corresponding results of [5] and [9] are immediate corollaries. We also define the rate function directly in terms of the $L_2(m)$-generator when $T = [0, \infty)$. If the semigroup is symmetric, so that $m$ is a $\sigma$-finite invariant measure, we identify the rate function in terms of the associated Dirichlet form. Our proof has the same outline as in [3], but we hope it is more transparent. We prove a lemma (Proposition 4.7) for the lower bound which would simplify the proof in the setting of [6] when the process is $m$-irreducible; we also do away with the need to define stopping times, as was done there.

We will always assume $E$ to be a topological Lusin space, which means that $E$ is a Hausdorff space and there exists a Polish (complete separable metric) space $S$, together with a one-to-one onto continuous map $\varphi: S \to E$. Let $\mathcal{B}(S)$ and $\mathcal{B}(E)$ denote the $\sigma$-algebras generated by open subsets of $S$ and $E$, respectively. It is known that $\varphi^{-1}$ is measurable (cf. [1]), therefore the $\sigma$-algebras $\mathcal{B}(S)$ and $\mathcal{B}(E)$ are $\sigma$-isomorphic. We also assume that $\mathcal{B}(E)$ is generated by $C(E)$, the real-valued continuous functions on $E$.

We will use the notation that if $\mathcal{A}$ denotes a class of (extended) real-valued measurable functions on some measurable space, then $\mathcal{A}^+$ denotes the nonnegative members of $\mathcal{A}$ and $\mathcal{A}_b$ denotes the bounded ones. We will also use $1_A(x) = 1(A)(x) = 1(0)$ if $x \in A (\notin A)$.

In addition to (1.1), we assume that the semigroup is generated by a transition probability function $p_t(x, \cdot)$, where for fixed $t \geq 0$, $x \in E$, $p_t(x, \cdot)$ is a probability measure on $(E, \mathcal{B}(E))$ and for $A \in \mathcal{B}(E)$, $p_t(x, A)$ is a jointly measurable function of $(t, x)$ with respect to $\mathcal{B}([0, \infty)) \times \mathcal{B}(E)$, where $\mathcal{B}([0, \infty))$ denotes the Borel subsets of $[0, \infty)$; if $t = 0$, then $p_t(x, A) = 1_A(x)$. For $f \in B_b(E)$, we write

$$p_t f(x) = \int f(y) p_t(x, dy)$$

and the semigroup $\{T_t, t \geq 0\}$ is then given on $L_2(m)$ via

(1.2) $$T_t f(x) = p_t f(x), \qquad f \in \hat{L}_2(m) \cap B_b(E),$$

in the sense that the function on the right-hand side is in the equivalence class that the left-hand side represents.

The $[L_2(m)]$ infinitesimal generator $L$ of the semigroup $\{T_t, t \in [0, \infty)\}$ is defined in the usual manner: $\mathcal{D}(L)$, the domain of $L$, consists of those members $f \in L_2(m)$ for which

$$\lim_{t \to 0} \frac{T_t f - f}{t} = g$$

exists in the $L_2(m)$-norm, and then $Lf := g$. In view of (1.1), $\mathcal{D}(L)$ is dense in $L_2(m)$. We define a resolvent function

(1.3) $$Rf(x) = \int_0^\infty e^{-t} p_t f(x) \, dt, \qquad f \in B_b(E),$$

and note an important consequence of (1.1) and (1.2), which is that $f = 0$ $m$-a.e. (almost everywhere with respect to $m$) implies that $Rf = 0$ $m$-a.e.

Next, we define an $m$-thin set.

DEFINITION 1.4. A set $N \in \mathcal{B}(E)$ is called $m$-thin if there exists a set $B \in \mathcal{B}(E)$ such that $m(B) = 0$ and

$$N \subset \{x \in E : R1_B(x) > 0\}.$$

Clearly, an $m$-thin set is $m$-null and it is also easy to see that a countable union of $m$-thin sets is $m$-thin.

We also need to make assumptions concerning the associated Markov process. When dealing with discrete time, we will take the path space $\Omega$ and the associated $\sigma$-algebras to have the product space representation. When $T = [0, \infty)$, we assume that there exists a progressively measurable Markov process $\{\Omega, \mathcal{F}_t^o, X_t, P^x, t \in T, x \in E\}$ with transition probability function $p_t(x, \cdot)$, where $X_t(\omega) = \omega(t)$ for $\omega \in \Omega$, $t \in T$, $\mathcal{F}_t^o = \sigma\{X_s, s \le t\}$ and $P^x(X_0 = x) = 1$ for each $x \in E$. The lifetime of the process starting from each $x$ is assumed to be $+\infty$. This structure will suffice for the local upper bound results. For the lower bound, we will use the entropy approach of Donsker and Varadhan [4] which requires $E$ to be Polish and $\Omega$ to be the corresponding Skorokhod space, which is also Polish. However, we would also like to include the setup of Ma and Röckner [7] and the results of [5] and [9] as special cases. For this reason, as far as the lower bound results are concerned, we assume that there exists a Polish space $\tilde{E}$ such that $E \in \mathcal{B}(\tilde{E})$, the Borel subsets of $\tilde{E}$. The space $\tilde{\Omega}$ is the Skorokhod space (the set of functions from $[0, \infty)$ to $\tilde{E}$, which are right-continuous and have left limits, with the Skorokhod topology). We assume that a Markov process $\{\tilde{\Omega}, \tilde{\mathcal{F}}_t^o, \tilde{X}_t, \tilde{P}^x, t \in T, x \in \tilde{E}\}$ exists with transition probability function $\tilde{p}_t(x, \cdot)$, where $\tilde{p}_t(x, A) = p_t(x, A \cap E)$ if $A \in \mathcal{B}(\tilde{E})$, $x \in E$, $t \ge 0$ and where $\tilde{p}_t(x, \{x\}) = 1$ if $x \notin E$, $t \ge 0$. Furthermore, we assume that there exists a Borel subset $\Omega$ of $\tilde{\Omega}$ such that if $\omega \in \Omega$, then $\tilde{X}_t(w) \in E$ for all $t \ge 0$ and if $x \in E$, then $P^x(\Omega) = 1$. We also assume that $P^x(\tilde{X}_0 = x) = 1$, and that the lifetime of the process is $+\infty$, for all $x \in \tilde{E}$. If $E$ itself is Polish, then the above assumption simply concerns the existence of the appropriate Markov process. As usual, $\theta_t$, $t \ge 0$, will denote the shift operator on the path space, that is, $\theta_t \omega(s) = \omega(s + t)$. For the case where, more generally, $E$ is a topological Lusin space, we will briefly describe the framework of Ma and Röckner [7] in the next paragraph.

In addition to (1.1) and (1.2), they assume in [7] that the semigroup is $m$-contractive, that is $\|T_t f\|_2 \le \|f\|_2$ for $f \in L_2(m)$. This assumption makes $-L$ a nonnegative definite operator. If we let

$$\mathcal{E}(u, v) = \langle -Lu, v \rangle, \qquad u, v \in \mathcal{D}(L),$$

then $\mathcal{E}$ is a nonnegative definite bilinear form. Define the norm

$$\|u\|_{\mathcal{E}}^2 = \mathcal{E}(u, u) + \|u\|_2^2, \qquad u \in \mathcal{D}(L),$$

and define $\mathcal{D}(\mathcal{E})$ to be the completion of $\mathcal{D}(L)$ under this norm. By the term "associated Dirichlet form" we mean the bilinear form $\mathcal{E}$ with domain $\mathcal{D}(\mathcal{E})$, which is required to satisfy the conditions of Definition 4.5, page 34 of [7]. The quasiregularity of $\mathcal{E}$ is defined on page 101 of [7]. For the sets of "capacity 0" or the notion of "quasi-everywhere," we again refer to Chapter III of [7]. Roughly speaking, $A \in \mathcal{B}(E)$ has capacity 0 if the process starting from $m$ as the initial measure never visits the set $A$. Under this framework, or where the notion of capacity in terms of the potential theory of the Markov process is available (which may happen even if $\mathcal{E}$ is not quasiregular), we will only need the following.

CONDITION 1.5. *If $f = 0$ $m$-a.e., then for $t > 0$, $p_t f = 0$, except on a set of capacity* 0.

This condition is satisfied under the quasiregularity condition on $\mathcal{E}$ imposed in [7] and implies that if a set is $m$-thin, then it has capacity 0 (cf. Lemma 2.19). If the semigroup is strongly Feller, that is, $T_t : B_b(E) \to C_b(E)$ for $t > 0$, and if for any $A \in \mathcal{B}(E)$, $m(A) = 0$ implies that the closure of $A^c$ (= complement of $A$) is $E$, then any $m$-thin set is empty. If $E = \mathbb{R}^n$, a large class of Lévy processes, many diffusion processes and ergodic stationary Markov processes have strongly Feller semigroups and satisfy the property mentioned above. Our results apply to these processes, even though their generators may not satisfy the sector condition (given by (2.5), page 16 of [7]) which is used for constructing the "right" Markov processes in [7].

To describe our results, we first define the empirical measures: if $A \in \mathcal{B}(E)$, $\omega \in \Omega$, then

$$L_t(\omega, A) := \begin{cases} \dfrac{1}{t} \displaystyle\int_0^t 1_A(X_s(\omega)) \, ds, & t \in (0, \infty), \\ \dfrac{1}{n} \displaystyle\sum_{j=0}^{n-1} 1_A(X_j(\omega)), & n \in \{1, 2, \dots\}. \end{cases}$$

The following subsets of $B(E)$ will play an important role:

$$D := \left\{ u \in B_b^+(E) : \int u \, dm < \infty \right\}.$$

Note that $D$ is dense in $L_1^+(m)$ and hence in $L_2^+(m)$.

$$D_0 := \left\{ v : v = t^{-1} \int_0^t p_s u \, ds, \text{ for some } u \in D, \text{ some } t > 0 \right\}.$$

$$D_1 := \{ p_\eta v : v \in D_0, \eta \geq 0 \}.$$

Members of $D, D_0$ and $D_1$ are functions, *not* equivalence classes. Also, note that $D_0 \subset D_1$.

We will see in Lemma 2.2, given in the next section that if $v \in D_1$ so that $v = t^{-1} \int_0^t p_s \tilde{u} \, ds$, where $\tilde{u} = p_\eta u$ for some $\eta \geq 0$ and $u \in D$, and $t > 0$, then $v \in \mathcal{D}(L)$ and

$$(1.6) \qquad \hat{L}v := \frac{1}{t}(p_t \tilde{u} - \tilde{u})$$

is an $m$-version of $Lv$.

We now define the rate function. $\mathcal{M}(E)$ denotes the set of probability measures on $(E, \mathcal{B}(E))$. For $\mu \in \mathcal{M}(E)$, we define

$$\hat{I}(\mu) := - \inf_{\substack{v \in D_1 \\ \varepsilon > 0}} \int \frac{\hat{L}v}{v + \varepsilon} \, d\mu \qquad \text{when } T = [0, \infty)$$

$$:= - \inf_{\substack{u \in B_b(E) \\ \varepsilon > 0}} \int \log \frac{p_1(u + \varepsilon)}{(u + \varepsilon)} \, d\mu \qquad \text{when } T = \{0, 1, 2, \dots\}$$

and

$$I(\mu) := \begin{cases} \hat{I}(\mu), & \text{if } \mu \ll m, \\ \infty, & \text{otherwise.} \end{cases}$$

A $\tau$-neighborhood of $\mu \in \mathcal{M}(E)$ is a set $N_\mu$ given by

(1.7) $$\left\{ \nu \in \mathcal{M}(E) : \left| \int f_j \, d\mu - \int f_j \, d\nu \right| < \varepsilon, 1 \le j \le r \right\}$$

for some $\varepsilon > 0$ and $f_1, \ldots, f_r$ in $B_b(E)$. This neighborhood system forms a basis for the $\tau$-topology on $\mathcal{M}(E)$. We denote by $\mathcal{B}_\tau(\mathcal{M}(E))$ the $\sigma$-algebra generated by these neighborhoods (*not* by all $\tau$-open sets). If the $f_j$'s are in $C_b(E)$, then $N_\mu$ is a weak neighborhood. A $w$-open subset of $\mathcal{M}(E)$ is one which is the union of weak neighborhoods of the form (1.7).

Before stating the results, we would like to note that the proofs for discrete time are obvious modifications of the ones that we are going to give for continuous time. In fact, there are fewer technicalities to deal with in discrete time. We will therefore simply state the results for discrete time in Section 6 without proofs.

We now state the local upper bound results when $T = [0, \infty)$, $\{T_t, t \ge 0\}$ satisfies (1.1) and (1.2), $E$ is a topological Lusin space and the Markov process is progressively measurable.

THEOREM 1.8 (Local upper bound). *Let $\mu \in \mathcal{M}(E)$ and let $a < \hat{I}(\mu)$. Then there exists a $\tau$-neighborhood $N_\mu$ of $\mu$ such that*

$$\limsup_{t \to \infty} \frac{1}{t} \sup_{x \in E} \log P^x(L_t \in N_\mu) \le -a.$$

We will see in Lemma 2.1 that $\tau$-compact sets are elements of $\mathcal{B}_\tau(M(E))$ and we obtain the following corollary of the theorem whose derivation uses standard arguments.

COROLLARY 1.9. *If $K$ is a $\tau$-compact subset of $\mathcal{M}(E)$, then*

$$\limsup_{t \to \infty} \frac{1}{t} \sup_{x \in E} \log P^x(L_t \in K) \le -\inf_{\mu \in K} \hat{I}(\mu).$$

REMARK 1.10. Note that $I(\mu) \ge \hat{I}(\mu)$ and that if $\mu$ is not absolutely continuous with respect to $m$, then $I(\mu) = +\infty$, in which case the lower bound in Theorem 1.15 below trivially holds and $\mu \ll m$ is the case of interest. In this sense, $I(\mu)$ is the natural rate function for the lower bound and the question naturally arises as to whether the upper bound holds in terms of $I(\mu)$. The next theorem answers this question.

THEOREM 1.11.   *Under the assumptions of Theorem 1.8, let* $\mu \in \mathcal{M}(E)$ *and* $a < I(\mu)$. *Then there exist a* $\tau$-*neighborhood* $N_\mu$ *of* $\mu$ *and an* $m$-*thin set* $N$ *such that if* $x \notin N$, *then*

$$\limsup_{t \to \infty} \frac{1}{t} \log P^x(L_t \in N_\mu) \leq -a.$$

*Furthermore, if* (1.5) *holds, then* capacity$(N) = 0$.

REMARK 1.12.   The $m$-thin set in the theorem depends on $N_\mu$ and one can give simple examples (cf. Example 3.5) to show that it need not be empty. The same example shows that there exists $\mu$, not absolutely continuous with respect to $m$, such that $\hat{I}(\mu) < \infty$, but the uniform upper bound does not hold if $\hat{I}(\mu)$ is replaced by $I(\mu)$.

REMARK 1.13.   The corollary of Theorem 1.11 corresponding to Corollary 1.9 also remains valid with an $m$-thin (or capacity 0) exceptional set.

We now state the lower bound results for which we have already described the setup. In addition to conditions (1.1) and (1.2), we will need the following.

CONDITION 1.14.   $m(A) > 0 \Rightarrow R\mathbf{1}_A(x) > 0, m$-a.e.$(x)$.

Even though we are considering $L_2(m)$-semigroups, we would like to observe that our rate functions are defined in terms of the transition probability function $p_t(x, \cdot)$. The irreducibility condition used in [6] says that there is a $\sigma$-finite reference measure $m$ on $(E, \mathcal{B}(E))$ such that if $A \in \mathcal{B}(E)$ and $m(A) > 0$, then (1.14) holds for *every* $x$. Under this stronger condition and an additional assumption that if $A$ is $m$-null, then $p_t\mathbf{1}_A(x)$ is zero for $m$-a.e. $(x)$, the exceptional set for large deviation lower bounds is shown to be empty (cf. other references in [6] as well). Under the weaker condition (1.14) together with the same additional assumption, we show here that the corresponding exceptional set is $m$-thin. The results of [6] follow from this result by simply observing that under the irreducibility assumption, starting from an arbitrary point $x$, the complement of the exceptional set will be visited at some time $t > 0$ with a positive probability and a routine argument then shows that the lower bound result must hold starting from *any* $x$. Our arguments here are simpler than those given in [6], in that we do not need to use stopping times and therefore the strong Markov property of the process is not needed.

THEOREM 1.15 (Lower bound). *Assume* (1.1), (1.2) *and* (1.14). *Let* $\mu \in \mathcal{M}(E)$ *and let* $N_\mu$ *be a $\tau$-neighborhood of $\mu$. Then there exists an $m$-thin set $N$ such that if $x \notin N$, then*

$$\liminf_{t \to \infty} \frac{1}{t} \log P^x(L_t \in N_\mu) \geq -I(\mu).$$

*Furthermore, if* (1.5) *holds, then* capacity$(N) = 0$.

By Lemma 2.1 in the next section, a $w$-open set $U$ belongs to $\in \mathcal{B}_\tau(\mathcal{M}(E))$ and we get the following corollary.

COROLLARY 1.16. *Let $U$ be a $w$-open subset of $\mathcal{M}(E)$. Then there exists an $m$-thin set $N$ [capacity$(N) = 0$, if* (1.5) *holds] such that if $x \notin N$, then*

$$\liminf_{t \to \infty} \frac{1}{t} \log P^x(L_t \in U) \geq - \inf_{\mu \in U} I(\mu).$$

If the semigroup is symmetric, in which case the $L_2$-generator $L$ is self-adjoint, then we obtain an appropriate form of Theorem 5 of [3] for the formula for $I(\mu)$ in terms of the associated symmetric Dirichlet form. In this case, $m$ is a $\sigma$-finite invariant measure for the semigroup.

THEOREM 1.17. *If $H = (-L)^{1/2}$ denotes the canonical square root of the positive definite operator $-L$, then for any $\mu \in \mathcal{M}(E)$ such that $\mu \ll m$, $I(\mu) < \infty$ if and only if $f = d\mu/dm \in \mathcal{D}(H)$, in which case*

$$I(\mu) = \|Hf\|_2^2.$$

REMARK 1.18 (Feller semigroups). If $\{T_t\}$ is a Feller semigroup, that is, $T_t : C_b(E) \to C_b(E)$ for each $t > 0$, in discrete or continuous time, in addition to being an $L_2(m)$ semigroup, Theorem 1.8 can be strengthened to the effect that $N_\mu$ may be chosen to be a $w$-neighborhood. As will be noted in the proof, no essential change in proof will be necessary if we modify the definition of the rate function $\hat{I}$ appropriately. Once $\hat{I}$ has been modified, $I$ is defined in terms of $\hat{I}$ as before. For $h > 0$, it will be useful to define $\hat{I}_h(\mu)$ by

$$(1.19) \qquad \hat{I}_h(\mu) = - \inf_{\substack{u \in C_b(E) \\ \varepsilon > 0}} \int \log \frac{p_h(u + \varepsilon)}{u + \varepsilon} \, d\mu$$

and note that (even if the semigroup is not Feller)

$$(1.20) \qquad \hat{I}_h(\mu) = - \inf_{\substack{u \in B_b^+(E) \\ \varepsilon > 0}} \int \log \frac{p_h(u + \varepsilon)}{u + \varepsilon} \, d\mu.$$

To see (1.20), given any $h > 0$, any $\mu \in \mathcal{M}(E)$ and any $\varepsilon > 0$, let

$$\mathcal{H} = \{u \in B_b(E) : \exists u_n \in C_b(E) \text{ such that } J(u_n; h, \varepsilon) \to J(u; h, \varepsilon)\},$$

where $J(u; h, \varepsilon)$ denotes the integral on the right-hand side of (1.19). It is easily verified that $\mathcal{H}$ satisfies the conditions of T20, page 11 of [8], and it follows that $\mathcal{H} = B_b(E)$ since $\mathcal{H} \supset C_b(E)$ and we are assuming that $\mathcal{B}(E)$ is generated by $C_b(E)$. In discrete time, we then define $\hat{I}(\mu) = \hat{I}_1(\mu)$ as before and no modification of the definition of the rate function is necessary; consequently, the upper (Theorem 1.8 and its Corollary 1.9) and the lower bound results will hold if "$\tau$-neighborhood" and "$\tau$-compact" are replaced by "$w$-neighborhood" and "$w$-compact," respectively.

In continuous time, the situation is slightly more complicated. This time, we will need to take

$$D = \left\{ u \in C_b^+(E) : \int u \, dm < \infty \right\}$$

and *assume* that $D$ is dense in $L_1^+(m)$ [in the $L_1(m)$-norm], hence also in $L_2^+(m)$. (Note that this assumption is satisfied if $E$ is Polish.) Then $D_0$ and $D_1$ are defined in terms of $D$ as before and the new rate function $\hat{I}(\mu)$ is defined as

$$\hat{I}(\mu) = - \inf_{\substack{v \in D_1 \\ \varepsilon > 0}} \int \frac{\hat{L}v}{v + \varepsilon} \, d\mu.$$

It will be shown in Lemma 2.15 that this new definition of $\hat{I}(\mu)$ agrees with our old definition [provided the modified $D$ is dense in $L_1^+(m)$]. We note that this requirement on $D$ is satisfied in many "usual" situations, where either $m$ is a probability measure or $E$ is a metric space and $m(K) < \infty$ for any compact $K$ and whenever open $U_n \searrow K$, $m(U_n) \searrow m(K)$. Once the new definition of $\hat{I}(\mu)$ agrees with the one given earlier, the proof of Theorem 1.8 needs an obvious modification, as noted there, so that $N_\mu$ may be chosen to be a $w$-neighborhood. Corollary 1.9 is valid for $w$-compact sets, and for the lower bound, Theorem 1.15 and its corollary remain valid *as stated*.

We would like to make some further comments about the methods of proof. As far as the upper bounds are concerned, we follow the approach of Donsker and Varadhan—use a form of the Feynman–Kac formula and apply the Chebyshev inequality. The form of this formula, which is very general, has been taken from Deuschel and Stroock [2]. Our lower bound result follows the method of Donsker and Varadhan [4], where they introduce the entropy of a stationary process relative to the given Markov process. A careful examination of the construction in [7] of a "right Markov process" with state space a topological Lusin space allows us to reduce this situation to the one where the state space and $\Omega$ are Polish, which is needed in [4]. We would like to note that dealing with the more general $E$ is not where the problem lies, one also needs $\Omega$ to be Polish to apply the results of [4] and the construction of appropriate $\Omega$ and Markov process requires considerable work [7], even when $E$ is Polish.

We present some preliminary results in Section 2. The upper bound results are in Section 3 and those for lower bounds are in Section 4. Section 5 deals with the symmetric situation (Theorem 1.17). Finally, the results for discrete time are stated in Section 6, without proofs.

## 2. Some preliminary results.

We present some useful lemmas in this section.

LEMMA 2.1. *Let $E$ be a topological Lusin space and let $\mathcal{B}_\tau(\mathcal{M}(E))$ and $\mathcal{B}_w(\mathcal{M}(E))$ denote the $\sigma$-algebras of subsets of $\mathcal{M}(E)$ generated by the $\tau$ and weak neighborhoods of the form (1.7), respectively. Then the $w$-open subsets of $\mathcal{M}(E)$ belong to $\mathcal{B}_\tau(\mathcal{M}(E))$. Consequently, $w$-closed and $\tau$-compact subsets of $\mathcal{M}(E)$ belong to $\mathcal{B}_\tau(\mathcal{M}(E))$.*

PROOF. Let $S$ be a Polish space and let $\varphi : S \to E$ be a continuous one-to-one and onto map. Let $\Phi(\mu) := \mu \varphi^{-1}$ if $\mu \in \mathcal{M}(S)$. It is then easy to check (note, again, that $\varphi^{-1}$ is measurable) that $\Phi : \mathcal{M}(S) \to \mathcal{M}(E)$ is a one-to-one and onto continuous map if $\mathcal{M}(S)$ and $\mathcal{M}(E)$ are given the $w$-topology. This makes $\mathcal{M}(E)$ a topological Lusin space and $\Phi^{-1}$ a measurable map, that is, $\Phi(\mathcal{B}_w(\mathcal{M}(S))) \subset \mathcal{B}_w(\mathcal{M}(E))$. If $U$ is a $w$-open subset of $\mathcal{M}(E)$, then $\Phi^{-1}(U)$ is $w$-open in $\mathcal{M}(S)$, which implies that $\Phi^{-1}(U) \in \mathcal{B}_w(\mathcal{M}(S))$ and its image $U$ under $\Phi$ belongs to $\mathcal{B}_w(\mathcal{M}(E)) \subset \mathcal{B}_\tau(\mathcal{M}(E))$. This proves the assertion about $w$-open sets. It follows that $w$-closed subsets of $\mathcal{M}(E)$ are in $\mathcal{B}_\tau(\mathcal{M}(E))$. Since a $\tau$-compact subset of $\mathcal{M}(E)$ is $w$-compact, it belongs to $\mathcal{B}_\tau(\mathcal{M}(E))$ and the lemma is proved. □

LEMMA 2.2. *Let $v \in D_1$ so that for some $u \in D$, some $t > 0$ and some $\eta \geq 0$, $v = t^{-1} \int_0^t p_s \tilde{u} \, ds$, where $\tilde{u} = p_\eta u$. Then for $h > 0$, $|h^{-1}(p_h v - v)| \leq c$ for some $c > 0$ independent of $h$ and*

$$\hat{L}v := t^{-1}(p_t \tilde{u} - \tilde{u})$$

*is bounded by $c$ and is an $m$-version of $Lv$.*

PROOF. We have

$$\frac{1}{h}(p_h v - v) = \frac{1}{t}\left\{\frac{1}{h}\left(\int_t^{t+h} p_s \tilde{u} \, ds - \int_0^h p_s \tilde{u} \, ds\right)\right\}.$$

Since $u$ is in $D$, it follows that $0 \leq u \leq c_1$ for some $c_1 > 0$ and the boundedness assertions follow. The strong continuity of the semigroup implies that as $h \to 0$, the right-hand side converges in $L_2(m)$ to $t^{-1}(p_t\tilde{u} - \tilde{u})$ and the lemma is proved. $\square$

LEMMA 2.3. $D_0$ is dense in $L_2^+(m)$. If the semigroup is symmetric, then $D_0$ is a dense subset of $L_1^+(m)$.

PROOF. Since $D$ is dense in $L_2^+(m)$, for the first assertion, it suffices to show that if $u \in D$, then there exist $v_n \in D_0$ such that $v_n \to u$ in $L_2(m)$. Let $t_n \downarrow 0$ and let $v_n = t_n^{-1}\int_0^{t_n} p_s u \, ds$. By the strong continuity of the semigroup, $\|v_n - u\|_2 \to 0$, so the first assertion is proved. For the second assertion, the symmetry of the semigroup implies that $m$ is an invariant measure for $\{T_t\}$, and for $u \in D$, $\int p_h u \, dm = \int u \, dm < \infty$. It follows that $D_0 \subset D$. Since $D$ is dense in $L_1^+(m)$, it suffices to show that if $u \in D$, then there exist $v_n \in D_0$ such that $\int |v_n - u| \, dm \to 0$ as $n \to \infty$. Taking $v_n$ in terms of $u$ as above, we have

$$\int |v_n - u| \, dm = 2 \int (u - v_n)^+ \, dm - \int (u - v_n) \, dm.$$

The second term on the right-hand side is 0 and $(u - v_n)^+ \leq u \in L_1^+(m)$ for all $n$; since $v_n \to u$ in $m$-measure by the strong continuity of the semigroup, $\int (v_n - u)^+ \, dm \to 0$, and the lemma is proved. $\square$

The next lemma is a version of the Feynman–Kac formula (see Deuschel and Stroock [2], page 121).

LEMMA 2.4. If $\varphi$ and $V$ are in $B_b(E)$, then the equation

$$(2.5) \qquad u(t,x) = p_t\varphi(x) + \int_0^t p_{t-s}(V(\cdot)u(s,\cdot))(x) \, ds$$

has a unique solution $u(t,x), t \geq 0, x \in E$, such that $\sup_{0 \leq t \leq a} \|u(t,x)\|_\infty < \infty$ for all $a \geq 0$.

The following corollary of this lemma is what we really need.

COROLLARY 2.6. Let $v \in D_1$ so that for some $\eta > 0$, $u \in D$ and $t_0 > 0$,

$$v = \frac{1}{t_0}\int_0^{t_0} p_s\tilde{u} \, ds,$$

where $\tilde{u} = p_\eta u$. For $\varepsilon > 0$, let

$$(2.7) \qquad v_\varepsilon(x) = v(x) + \varepsilon$$

and

$$(2.8) \qquad V_\varepsilon(x) := \frac{1}{v_\varepsilon(x)}\frac{p_{t_0}\tilde{u}(x) - \tilde{u}(x)}{t_0} = \frac{\hat{L}v}{v_\varepsilon}(x).$$

*Then for all $t > 0$, $x \in E$, we have*

$$(2.9) \qquad v_\varepsilon(x) = E^x\{v_\varepsilon(X_t)e^{-\int_0^t V_\varepsilon(X_s)\,ds}\}.$$

PROOF. Let $v_\varepsilon(t,x)$ denote the right-hand side of (2.9). By (4.2.25) in [2], we get that $v_\varepsilon(t,x)$ is the unique solution of

$$(2.10) \qquad v_\varepsilon(t,x) = p_t v_\varepsilon(x) - \int_0^t p_{t-s}(V_\varepsilon(\cdot)v_\varepsilon(s,\cdot))(x)\,ds.$$

We will show that $v_\varepsilon(t,x) = v_\varepsilon(x)$ satisfies (2.10). Substituting for $V_\varepsilon$ and $v_\varepsilon(s,x) = v_\varepsilon(x)$ from (2.7) in the right-hand side of (2.10), we get

$$p_t v_\varepsilon(x) - \int_0^t p_s\left(\frac{1}{t_0}(p_{t_0}\tilde{u} - \tilde{u})\right)(x)\,ds$$

and by the definition of $v$, this equals

$$\varepsilon + \frac{1}{t_0}\int_0^{t_0} p_{t+s}\tilde{u}(x)\,ds - \frac{1}{t_0}\int_0^t p_s(p_{t_0}\tilde{u} - \tilde{u})(x)\,ds = \varepsilon + v(x) = v_\varepsilon(x),$$

which proves the corollary.  □

LEMMA 2.11. *Let $\mu \in \mathcal{M}(E)$, $\mu \ll m$. Given $h > 0$, and writing $u_\varepsilon = u + \varepsilon$ for $\varepsilon > 0$, $u \in B_b^+(E)$, let*

$$\theta(u;h,\varepsilon) := \log(p_h u_\varepsilon / u_\varepsilon), \qquad J_\mu(u;h,\varepsilon) := \int \theta(u;h,\varepsilon)\,d\mu.$$

*We then have (Feller case included)*

$$(2.12) \qquad \inf_{\substack{v \in D_0 \\ \varepsilon > 0}} J_\mu(v;h,\varepsilon) = \inf_{\substack{v \in D_1 \\ \varepsilon > 0}} J_\mu(v;h,\varepsilon) = \inf_{\substack{u \in B_b^+(E) \\ \varepsilon > 0}} J_\mu(u;h,\varepsilon).$$

PROOF. Note that $D_0 \subset D_1 \subset B_b^+(E)$, so it suffices to show that the first and last terms are equal. We first consider the general case without assuming that the semigroup is Feller. The left-hand side in (2.12) is larger than the right-hand side. Let $u \in B_b^+(E)$. Since $m$ is $\sigma$-finite, there exist $A_n \nearrow E$ such that $m(A_n) < \infty$ for $n \geq 1$. Let $u_n = u1_{A_n}$. Then $u_n \in D$. Clearly, $\theta(u_n;h,\varepsilon) \to \theta(u;h,\varepsilon)$ boundedly as $n \to \infty$ for a fixed $\varepsilon > 0$. Therefore, $J_\mu(u_n;h,\varepsilon) \to J_\mu(u;h,\varepsilon)$ as $n \to \infty$. It follows that (2.12) holds if $D_0$ is replaced by $D$. To complete the proof, it suffices to show that

$$(2.13) \qquad \inf_{\substack{v \in D_0 \\ \varepsilon > 0}} J(v;h,\varepsilon) \leq \inf_{\substack{u \in D \\ \varepsilon > 0}} J(u;h,\varepsilon).$$

To see this, let $u \in D$ and for some $t_n \downarrow 0$, let $v_n = t_n^{-1}\int_0^{t_n} p_s u\,ds$. Then $v_n \in D_0$ and by the strong continuity of the semigroup in $L_2(m)$, $v_n \to u$ in $L_2(m)$. Since $v_n$ and $u$ are bounded by some $c > 0$ and $\mu \ll m$, we have $v_n \to u$ in $L_2(\mu)$, hence in $\mu$-measure, so $J_\mu(v_n; h, \varepsilon) \to J_\mu(u; h, \varepsilon)$ as $n \to \infty$ and the lemma is proved for the general case. We now consider the Feller case.

Let $D_0$ and $D$ be as in the general case and, for the moment, let $D_0' = \{v : v = \frac{1}{t} \int_0^t p_s u \, du$ for some $t > 0$ and some $u \in D \cap C_b(E)\}$. It then suffices to show that

$$\inf_{\substack{v \in D_0' \\ \varepsilon > 0}} J(v; h, \varepsilon) \leq \inf_{\substack{v \in D_0 \\ \varepsilon > 0}} J(v; h, \varepsilon). \tag{2.14}$$

Let $v \in D_0$. Then $v = \frac{1}{t} \int_0^t p_s u \, du$ for some $t > 0$ and some $u \in D$. By our assumption in the Feller case namely, that $D \cap C_b(E)$ is dense in $L_1^+$, it is easy to see ($\mu \ll m$) that there exist $u_n \in D \cap C_b(E)$ such that $u_n \to u$ boundedly in $\mu$-measure. Therefore, $v_n := \frac{1}{t} \int_0^t p_s u_n \, ds \to v$ boundedly in $\mu$-measure, hence $J(v_n; h, \varepsilon) \to J(v; h, \varepsilon)$ as $n \to \infty$. Since $v_n \in D_0'$, (2.14) is proved and (2.12) holds in the Feller case as well. $\square$

LEMMA 2.15.   *Let $\mu \in \mathcal{M}(E), \mu \ll m$. Then (including the Feller case)*

$$\lim_{h \to 0} \frac{1}{h} \hat{I}_h(\mu) = \hat{I}(\mu). \tag{2.16}$$

PROOF.   If $v \in D_1$, then there exists $c > 0$ such that for all $h > 0$,

$$\frac{1}{h} |p_h v - v| \leq c.$$

Therefore, writing $v_\varepsilon = v + \varepsilon$, we have

$$\log p_h v_\varepsilon = \log v_\varepsilon + (p_h v_\varepsilon - v_\varepsilon) \cdot \frac{1}{v_\varepsilon} + O(h^2),$$

where $O$ depends only on $\varepsilon$. As $h \to 0$, $h^{-1}(p_h v - v) \to \hat{L}v$ in $L_2(m)$ boundedly by Lemma 2.2, hence in $\mu$-measure boundedly, and we get

$$\limsup_{h \to 0} \frac{1}{h} \inf_{\substack{v \in D_1 \\ \varepsilon > 0}} \int \log \frac{p_h v_\varepsilon}{v_\varepsilon} \, d\mu \leq \int \frac{\hat{L}v}{v_\varepsilon} \, d\mu$$

for any $v \in D_1$, $\varepsilon > 0$. By Lemma 2.11, we get

$$\liminf_{h \to 0} \frac{1}{h} \hat{I}_h(\mu) \geq \hat{I}(\mu). \tag{2.17}$$

We now consider the opposite inequality and use an argument of Donsker and Varadhan. For $v \in D_0$, define

$$\varphi(h) = \int \log \frac{p_h v_\varepsilon}{v_\varepsilon} \, d\mu,$$

where $v_\varepsilon = v + \varepsilon$. Then $v \in \mathcal{D}(L)$ and we have

$$\frac{d\varphi}{dh} = \int \frac{\hat{L} p_h v}{p_h v_\varepsilon} \, d\mu \geq \inf_{\substack{v \in D_1 \\ \varepsilon > 0}} \int \frac{\hat{L} v}{v_\varepsilon} \, d\mu = -\hat{I}(\mu).$$

Integrating in $h$ from 0 to $h$ and using $\varphi(0) = 0$, we get (for all $v \in D_1$, $\varepsilon > 0$, $h > 0$)

$$(2.18) \qquad \int \log \frac{p_h v_\varepsilon}{v_\varepsilon} \, d\mu \geq -h \hat{I}(\mu),$$

which, together with (2.17) and Lemma 2.11, finishes the proof of the lemma. $\square$

When the potential-theoretic framework is available, such as in [7], the following lemma shows that an $m$-thin set has capacity 0.

LEMMA 2.19. *If* (1.5) *and the related framework hold, then an $m$-thin set $N$ has capacity 0.*

PROOF. Let $N = \{x \in E : R1_A(x) > 0\}$, where $m(A) = 0$. Then $u(x) = R1_A(x) = 0$ $m$-a.e. Hence, $p_s u(x) = \int_s^\infty e^{-t} p_t(x, A) \, dt = 0$, except on a set $N_s$ of capacity 0, by (1.5). Since $N = \bigcup_{j=1}^\infty N_{1/j}$ and capacity$(N_{1/j}) = 0$, we have capacity$(N) = 0$ because a countable union of sets of capacity 0 has capacity 0. The lemma is proved. $\square$

## 3. The local upper bound proofs.
In this section, we prove Theorems 1.8 and 1.11.

PROOF OF THEOREM 1.8. Let $v \in D_1$ and for $\varepsilon > 0$ let $v_\varepsilon = v + \varepsilon$. By Corollary 2.6, for $x \in E$, $t > 0$ and $\varepsilon > 0$, we have

$$v_\varepsilon(x) = E^x \left\{ v_\varepsilon(X_t) \exp\left( -\int_0^t \frac{\hat{L} v}{v_\varepsilon}(X_s) \, ds \right) \right\}.$$

Writing $V_\varepsilon := \hat{L} v / v_\varepsilon$, we have, for $A \in \mathcal{B}_\tau(\mathcal{M}(E))$,

$$
\begin{aligned}
v_\varepsilon(x) &= E^x \left\{ v_\varepsilon(X_t) \exp\left( -t \int V_\varepsilon \, dL_t \right) \right\} \\
&\geq \varepsilon E^x \left\{ \exp\left( -t \int V_\varepsilon \, dL_t \right) 1_{[L_t \in A]} \right\} \\
&\geq \varepsilon E^x \left\{ \exp\left( -t \sup_{\nu \in A} \int V_\varepsilon \, d\nu \right) 1_{[L_t \in A]} \right\}.
\end{aligned}
$$

Since $v_\varepsilon \le \alpha$ for some $0 < \alpha < \infty$, where $\alpha$ may depend on $\varepsilon$, we get

$$\frac{1}{t}\log\sup_x P^x(L_t \in A) \le \frac{\log\alpha - \log\varepsilon}{t} + \sup_{\nu \in A}\int V_\varepsilon\,d\nu.$$

Therefore, for any $\varepsilon > 0$ and any $v \in D_1$ any $A \in \mathcal{B}_\tau(\mathcal{M}(E))$, we have

$$(3.1) \qquad \limsup_{t\to\infty}\frac{1}{t}\log\sup_{x\in E}P^x(L_t \in A) \le \sup_{\nu\in A}\int\frac{\hat{L}v}{v+\varepsilon}\,d\nu.$$

Since $a < \hat{I}(\mu)$, for some $\varepsilon_0 > 0$ and some $v_0 \in D_1$, we have, for some $\delta > 0$,

$$\int\frac{\hat{L}v_0}{v_0+\varepsilon_0}\,d\mu = -a - \delta.$$

We now select $N_\mu = \{\nu \in \mathcal{M}(E): |\int f\,d\nu - \int f\,d\mu| < \frac{\delta}{2}\}$, where $f = \hat{L}v_0/(v_0 + \varepsilon_0)$, so $N_\mu$ is a $\tau$-neighborhood of $\mu$. In the Feller case, since the members of $D_1$ are in $C_b(E)$, $\hat{L}v_0 \in C_b(E)$. Consequently, this choice of $N_\mu$ automatically yields a $w$-neighborhood. In either case, we have

$$(3.2) \qquad \sup_{\nu\in N_\mu}\int\frac{\hat{L}v_0}{v_0+\varepsilon_0}\,d\nu < -a.$$

If we take $A = N_\mu$, $v = v_0$ and $\varepsilon = \varepsilon_0$ in (3.1) and use (3.2), then the theorem follows. $\square$

PROOF OF THEOREM 1.11. If $\mu \ll m$, then $I(\mu) = \hat{I}(\mu)$ and we get the result for all $x$, as in Theorem 1.8. Now, let $\mu$ have a nonzero singular component with respect to $m$. There then exists a set $A$ such that $\mu(A) > 0$ and $m(A) = 0$. In this case, $I(\mu) = +\infty$ and we will first show that there exist a $\tau$-neighborhood $N_\mu$ of $\mu$ and an $m$-thin set $N$ such that if $x \notin N$, then

$$(3.3) \qquad \limsup_{t\to\infty}\frac{1}{t}\log P^x(L_t \in N_\mu) = -\infty.$$

Let $\mu(A) = \eta > 0$ and define

$$N_\mu = \left\{\nu \in \mathcal{M}(E): \left|\int 1_A\,d\nu - \int 1_A\,d\mu\right| < \eta/2\right\}.$$

By (1.2), if $m(A) = 0$, then for all $t > 0$, $m$-a.e.$(x)$, $p_t(x, A) = 0$, hence

$$\int_0^\infty\int m(dx)p_s(x, A)\,ds = 0,$$

which implies that

$$\int m(dx)E^x\int_0^\infty 1_A(X_s)\,ds = 0.$$

If we let

$$N = \left\{ x : E^x \int_0^\infty 1_A(X_s)\, ds > 0 \right\},$$

then $N$ is $m$-thin. If $x \notin N$, then $\int_0^\infty p_t(x, A)\, dt = 0$ and we have, for all $t > 0$,

(3.4)                          $P^x(|L_t(A) - \mu(A)| < \eta/2) = 0,$

which proves (3.3). If Condition 1.5 holds, then by Lemma 2.19, the set $N$ which is $m$-thin must have capacity 0 and the theorem is proved.  □

EXAMPLE 3.5.   This example is a trivial modification of standard Brownian motion on $\mathbb{R}$, but shows the necessity of an exceptional set in Theorem 1.11 if $I(\mu)$ is to be used as the rate function. This example also shows that $\hat{I}(\mu) < \infty$ and, of course, Theorem 1.8 holds, but becomes false if $I(\mu)$ is used in place of $\hat{I}(\mu)$. Let $E = \mathbb{R} \cup \{\theta\}$, $\mathbb{R}$ have the usual topology of the real line and $\theta$ be an isolated point. Let $\{X_t, t \geq 0\}$ be the standard Brownian motion on $\mathbb{R}$. If $X_0 = x \in \mathbb{R}$, then so that $p_t(x, \mathbb{R}) = 1$ for $t \geq 0$ and if $X_0 = \theta$, then $p_t(\theta, \{\theta\}) = 1$. The Lebesgue measure on $\mathbb{R}$ plays the role of $m$ with $m(\{\theta\}) = 0$. Let $\mu = \delta_\theta$ and let $N_\mu = \{\nu : |\int f_j\, d\nu - \int f_j\, d\mu| < \varepsilon, 1 \leq j \leq r\}$ be any $\tau$-neighborhood of $\mu$. Then $p^\theta(L_t \in N_\mu) = 1$ and since $I(\mu) = +\infty$, the upper bound can hold only $m$-a.e. The exceptional set in this case is $\{\theta\}$, which is $m$-thin. Note also that $\hat{I}(\mu) = 0$, so the conclusion of Theorem 1.8, though not very interesting, holds.

**4. Proof of the lower bound.**   Recall that Condition 1.14 is assumed in addition to (1.1) and (1.2) for the lower bound.

If $E$ is Polish, then we assume that there exists a Markov process (as described in Section 1) with path space $\Omega$, the Skorokhod space corresponding to $E$. If $E$ is Lusin, then the Borel structure of $E$ is the same as that of a Polish space $S$ (see Section 1), but the topology of $E$ may not be metrizable. Under the main assumption that the associated Dirichlet form is quasiregular, Ma and Röckner [7] construct a "right" Markov process. In their construction, they imbed $E$ as a Borel subset into a Polish space $\tilde{E}$ in such a way that the "right" Markov process has $\tilde{\Omega}$ (the Skorokhod space corresponding to $\tilde{E}$) as the path space. Furthermore, there exists a Borel subset $\Omega$ of $\tilde{\Omega}$ such that the paths starting from $x \in E$ all lie in $\Omega$. To obtain the lower bound results and use the framework of [4], we simply work with the Markov process with state space $\tilde{E}$ and path space $\tilde{\Omega}$ and relativize these results to $E$ and $\Omega$. Therefore, in the setting of [7], there is no loss of generality in assuming $E$ and $\Omega$ to be Polish and we will do so.

We start with the following, simple, lemma which is a consequence of Condition 1.14.

LEMMA 4.1.   *Let $\nu \in \mathcal{M}(E)$, $\nu \ll m$. If $m(A) > 0$, then there exists a $t_0 > 0$ and a set $B \in \mathcal{B}(E)$ such that $\nu(B) > 0$ and*

$$\inf_{x \in B} p_{t_0}(x, A) > 0.$$

PROOF.   An immediate consequence of Condition 1.14 is that

$$\int \nu(dx) \int_0^\infty p_t(x, A)\, dt > 0,$$

therefore for some $t_0 > 0$, $\int p_{t_0}(x, A)\nu(dx) > 0$. The conclusion of the lemma immediately follows from this.   □

To establish the lower bound, we will need a basic result of Donsker and Varadhan [4] called the "contraction principle" and an approximation lemma from [6]. Let $\mathcal{M}_S(\Omega)$ denote the set of stationary probability measures on $(\Omega, \mathcal{F}^o)$. We state these results here for ready reference, but first, we need to define the entropy of $Q \in \mathcal{M}_S(\Omega)$ on $(\Omega, \mathcal{F}_t^o)$ with respect to the Markov process, as introduced in [4]. For $t > 0$, let

$$H(t, Q) := E^Q\{h(P^{\omega(0)}|\mathcal{F}_t^o; Q_\omega|\mathcal{F}_t^o)\},$$

where $h(\lambda; \mu)$ denotes the entropy of $\mu$ with respect to $\lambda$, $\lambda \,|\, \mathcal{G}$ denotes the restriction of $\lambda$ to $\mathcal{G}$, $Q_\omega$ is the regular conditional probability distribution of $Q$ with respect to $\mathcal{F}_o^o$ and $E^Q$ denotes the expectation with respect to the stationary measure $Q$ on $(\Omega, \mathcal{F}^o)$. The following facts about $H(t, Q)$ are established in [4]:

(i) either $H(t, Q) = +\infty$ for all $t > 0$ or there exists a nonnegative constant $H(Q)$ such that $H(t, Q) = tH(Q)$, and if $H(Q) < \infty$, then $Q$-a.s.($\omega$), $Q_\omega \ll P^{\omega(0)}$;

(ii) if $\lambda_1 \geq 0$, $\lambda_2 \geq 0$, $\lambda_1 + \lambda_2 = 1$ and $Q_1, Q_2 \in \mathcal{M}_S(\Omega)$, then $H(\lambda_1 Q_1 + \lambda_2 Q_2) = \lambda_1 H(Q_1) + \lambda_2 H(Q_2)$;

(iii) if $Q \in \mathcal{M}_S(\Omega)$ is ergodic and $\psi_t(\omega, \cdot)$ denotes the Radon–Nikodym derivative of $Q_\omega$ with respect to $P^{\omega(0)}$, both restricted to $\mathcal{F}_t^o$, then $H(Q) < \infty$ implies that $Q$-a.s.($\omega$), $Q_\omega$-a.s.,

(4.2)                    $$\lim_{t \to \infty} \frac{1}{t} \log \psi_t(\omega, \cdot) = H(Q).$$

THEOREM 4.3 (Contraction principle).   *Let $\mu \in \mathcal{M}(E)$, $\mu \ll m$ and $I(\mu) < \infty$. Then $I(\mu) = \inf\{H(Q) : Q \text{ stationary with marginal } \mu\}$, where $I(\mu)$ is defined in Section 1.*

PROOF.   It is proved in [4] that for any $\mu \in \mathcal{M}(E)$, $\inf\{H(Q) : Q \text{ stationary with marginal } \mu\} = \bar{I}(\mu)$, where

$$\bar{I}(\mu) := \limsup_{h \downarrow 0} \frac{1}{h} \hat{I}_h(\mu)$$

and $\hat{I}_h(\mu)$ is given by (1.19). If $\mu \ll m$, then by Lemma 2.15, we have $\bar{I}(\mu) = I(\mu)$. □

The next lemma is Lemma 2.5 of [6].

LEMMA 4.4.   Let $Q_0 \in \mathcal{M}_S(\Omega)$. Let $G$ be a $\tau$-neighborhood of $Q_0$ in $\mathcal{M}_S(\Omega)$,

$$G = \left\{ Q \in \mathcal{M}_S(\Omega) : \left| \int g_j \, dQ - \int g_j \, dQ_0 \right| < \varepsilon, 1 \leq j \leq r \right\},$$

where the $g_j$'s are bounded measurable functions on $(\Omega, \mathcal{F}^o)$. Then, given $\delta > 0$, there exists a $Q(\delta) = \sum_{p=1}^k \lambda_p Q_p$, where $\lambda_p > 0$, $Q_p \in \mathcal{M}_S(\Omega)$ is ergodic, $1 \leq p \leq k$, and $\sum_{p=1}^k \lambda_p = 1$, such that $Q(\delta) \in G$ and $|H(Q(\delta)) - H(Q_0)| < \delta$.

PROOF.   The proof in [6] is given when $G$ is a *weak* neighborhood. However, that proof is ergodic-theoretic and uses only the measurability and boundedness of the $f_j$'s, hence it applies without modification to our present situation. □

PROOF OF THEOREM 1.15.   Proposition 4.7 below is the first step in the proof. We write

$$N_{\mu,\varepsilon} = \left\{ \nu \in \mathcal{M}(E) : \left| \int f_j \, d\mu - \int f_j \, d\nu \right| < \varepsilon, 1 \leq j \leq r \right\}$$

for $\mu \in \mathcal{M}(E)$, $\varepsilon > 0$ and $f_1, \ldots, f_r \in B_b(E)$.

Note that it suffices to consider the case when $I(\mu) < \infty$ so that $\mu \ll m$. By the contraction principle, we can restrict our attention to those $Q_\mu$ for which $H(Q_\mu) < \infty$.

Before proceeding with the proof, we would like to make a remark which will be used in the proof more than once.

REMARK 4.5.   (i) Let

(4.6)        $$N_{\mu,\varepsilon} = \left\{ \nu \in \mathcal{M}(E) : \left| \int f_j \, d\nu - \int f_j \, d\mu \right| < \varepsilon, 1 \leq j \leq r \right\}$$

be a ($\tau$- or $w$-) neighborhood of $\mu$. Let $0 < \varepsilon_1 < \varepsilon$ and let $t_0 > 0$ be given. Then for all $t$ sufficiently large (depending only on $\varepsilon_1$ and $t_0$) and for all $\omega \in \Omega$, $L_t(\theta_{t_0}\omega) \in N_{\mu,\varepsilon_1}$ implies that $L_t(\omega) \in N_{\mu,\varepsilon}$. To see this, note that $L_t(\theta_{t_0}\omega) \in N_{\mu,\varepsilon_1}$ means

$$\left| \frac{1}{t} \int_{t_0}^{t_0+t} f_j(X_s(\omega)) \, ds - \int f_j \, d\mu \right| < \varepsilon_1, \qquad 1 \leq j \leq r,$$

which implies that

$$\left| \frac{1}{t} \int_0^t f_j(X_s(\omega)) \, ds - \int f_j \, d\mu \right| < \varepsilon_1 + \frac{2ct_0}{t}, \qquad 1 \le j \le k,$$

where $c = \max_{1 \le j \le r} \|f_j\|_\infty$, which proves the assertion.

(ii) The function

$$\varphi(x) := \liminf_{t \to \infty} \log E^x \left( 1(L_t \in N_{\mu,\eta}) \frac{1}{\varepsilon t} \int_t^{t(1+\varepsilon)} 1_A(X_s) \, ds \right),$$

where $N_{\mu,\eta}$ is a $w$- or a $\tau$-neighborhood of $\mu$ and $A \in \mathcal{B}(E)$, is a measurable function. The reason for this is that the lim inf can be taken along positive integers. Indeed, if $n \le t < n+1$, then for $1 \le j \le r$,

$$\frac{1}{n} \int_0^n f_j(X_s) \, ds - \frac{2c}{n} \le \frac{1}{t} \int_0^t f_j(X_s) \, ds \le \frac{1}{n} \int_0^n f_j(X_s) \, ds + \frac{2c}{n},$$

where $c$ is as in (i) above. Then, for all $t$ sufficiently large, $L_t \in N_{\mu,\eta}$ if and only if $L_n \in N_{\mu,\eta}$. Also,

$$\frac{1}{\varepsilon n} \int_n^{n(1+\varepsilon)} 1_A(X_s) \, ds - \frac{1+\varepsilon}{n} \le \frac{1}{\varepsilon t} \int_t^{t(1+\varepsilon)} 1_A(X_s) \, ds$$
$$\le \frac{1}{\varepsilon n} \int_n^{n(1+\varepsilon)} 1_A(X_s) \, ds + \frac{1+\varepsilon}{n},$$

which shows that the lim inf in the definition of $\varphi(x)$ can be taken along positive integers, hence $\varphi$ is measurable.

PROPOSITION 4.7.   *Let $\mu \in \mathcal{M}(E)$, $\mu \ll m$. Let $N_{\mu,\eta}$ be a $\tau$-neighborhood of $\mu$. Assume that $Q_\mu \in \mathcal{M}_S(\Omega)$ with marginal $\mu$ and that $Q_\mu$ is ergodic with $H(Q_\mu) < \infty$. If we let $A$ be such that $\mu(A) > 0$, then for $m$-a.e.$(x)$, we have*

$$(4.8) \quad \liminf_{t \to \infty} \frac{1}{t} \log E^x \left( 1(L_t \in N_{\mu,\eta}) \frac{1}{\varepsilon t} \int_t^{t(1+\varepsilon)} 1_A(X_s) \, ds \right) \ge -(1+\varepsilon) H(Q_\mu).$$

PROOF.   For any $\omega \in \Omega$,

$$(4.9) \quad \int_{(L_t \in N_{\mu,\eta})} \frac{1}{\varepsilon t} \int_t^{t(1+\varepsilon)} 1_A(X_s) \, ds \, dP^{\omega(0)}$$
$$\ge \int_{(L_t \in N_{\mu,\eta})} \frac{1}{\varepsilon t} \int_t^{t(1+\varepsilon)} 1_A(X_s) \, ds \frac{dP^{\omega(0)}}{dQ_{\mu,\omega}} \Big|_{\mathcal{F}_{t(1+\varepsilon)}^o} dQ_{\mu,\omega}.$$

[The Radon–Nikodym derivative here is that of the absolutely continuous part of $P^{\omega(0)}$ with respect to $Q_{\mu,\omega}$. We observed before that $Q_{\mu,\omega} \ll P^{\omega(0)}$, $Q_\mu$-a.s., if $H(Q_\mu) < \infty$, but the converse need not be true, hence the inequality.] Now, $Q_\mu$-a.s.($\omega$), as $t \to \infty$, by the ergodic theorem, we have

$$Q_{\mu,\omega}(L_t \in N_{\mu,\varepsilon}) \to 1,$$

$$\frac{1}{t(1+\varepsilon)} \log \frac{dP^{\omega(0)}}{dQ_{\mu,\omega}}\Big|_{\mathcal{F}^o_{t(1+\varepsilon)}} \to -H(Q_\mu)$$

and

$$\frac{1}{\varepsilon t} \int_t^{t(1+\varepsilon)} 1_A(X_s)\, ds \to \mu(A).$$

It follows that given $\delta > 0$, for all $t$ sufficiently large (depending on $\omega$ and $\delta$) $Q_\mu$-a.s.($\omega$), $Q_{\mu,\omega}$-a.s., the left-hand side of (4.9) is

$$\geq \frac{\mu(A)}{2}(1-\delta)\exp(-H(Q_\mu)-\delta)(t(1+\varepsilon)).$$

Therefore, for $\mu$-a.e.($x$), the left-hand side of (4.8) is

$$\geq -(1+\varepsilon)(H(Q_\mu)+\delta),$$

where the null set depends on $\delta$. Taking $\delta_n \downarrow 0$ and combining the countable number of null sets, we conclude that (4.8) holds. We now go from $\mu$-a.e. to $m$-a.e.

Since (4.8) holds $\mu$-a.e.($x$), given $\delta > 0$, there exist a $t_0 > 0$ and a set $B$ with $\mu(B) > 0$ such that for all $x \in B$, $t \geq t_0$, we have

$$(4.10)\quad \frac{1}{t} \log E^x\bigg(1(L_t \in N_{\mu,\eta})\frac{1}{\varepsilon t}\int_t^{t(1+\varepsilon)} 1_A(X_s)\, ds\bigg) \geq -(1+\varepsilon)H(Q_\mu)-\delta.$$

Since $\mu \ll m$, $\mu(B) > 0$ implies $m(B) > 0$. By Condition 1.14, we have $R1_B(x) > 0$ $m$-a.e.($x$). Let

$$C = \{x : R1_B(x) > 0\}.$$

If $x_0 \in C$, then there exists $t_1 > 0$ such that $p_{t_1}(x_0, B) = \alpha > 0$. For any $t > 0$, we have

$$\int_B p_{t_1}(x_0, dy) E^y\bigg(1(L_t \in N_{\mu,\eta})\frac{1}{\varepsilon t}\int_t^{t(1+\varepsilon)} 1_A(X_s)\, ds\bigg)$$

$$\leq E^{x_0}\bigg(1(L_t(\theta_{t_1}\cdot) \in N_{\mu,\eta})\frac{1}{\varepsilon t}\int_t^{t(1+\varepsilon)} 1_A(X_{s+t_1})\, ds\bigg)$$

and by Remark 4.5(i), for all sufficiently large $t$,

$$\leq E^{x_0}\bigg(1(L_t \in N_{\mu,2\eta})\frac{1}{\varepsilon t}\int_{t+t_1}^{(t+t_1)(1+\varepsilon)} 1_A(X_s)\, ds\bigg)$$

$$\leq E^{x_0}\bigg(1(L_t \in N_{\mu,2\eta})\frac{1}{\varepsilon t}\int_t^{t(1+\varepsilon)} 1_A(X_s)\, ds + \frac{t_1(1+\varepsilon)}{\varepsilon t}\bigg).$$

Taking logarithms and applying Jensen's inequality, we get

$$\frac{1}{\alpha} \int_B p_{t_1}(x_0, dy) \log E^y \Big( 1(L_t \in N_{\mu,\eta}) \frac{1}{\varepsilon t} \int_t^{t(1+\varepsilon)} 1_A(X_s) \, ds \Big)$$

$$\leq -\log \alpha + \log E^{x_0} \Big( 1(L_t \in N_{\mu,2\eta}) \frac{1}{\varepsilon t} \int_t^{t(1+\varepsilon)} 1_A(X_s) \, ds + \frac{t_1(1+\varepsilon)}{\varepsilon t} \Big).$$

If $y \in B$, and $t \geq t_0$, then the log term on the left-hand side divided by $t$ is bounded below [by (4.10)]. Hence, dividing by $t$, by Fatou's lemma, as $t \to \infty$, we get [via (4.10)] for $x_0 \in C$,

$$-(1+\varepsilon) H(Q_\mu) - \delta \leq \liminf_{t \to \infty} \frac{1}{t} \log E^{x_0} \Big( 1(L_t \in N_{\mu,2\eta}) \frac{1}{\varepsilon t} \int_t^{t(1+\varepsilon)} 1_A(X_s) \, ds \Big)$$

and since $\delta$ is arbitrary and the right-hand side does not depend on $\delta$, the proposition is proved. $\square$

The following corollary will be used in the next step.

COROLLARY 4.11. *Let* $\nu \in \mathcal{M}(E)$, $\nu \ll m$. *Then given* $\delta > 0$, *there exist a* $t_0 > 0$ *and a set* $B \in \mathcal{B}(E)$ *with* $\nu(B) > 0$ *such that for all* $t \geq t_0$,

$$\inf_{x \in B} E^x \Big( 1(L_t \in N_{\mu,\eta}) \frac{1}{\varepsilon t} \int_t^{t(1+\varepsilon)} 1_A(X_s) \, ds \Big) \geq \exp(-t(1+\varepsilon)(H(Q_\mu) + \delta)).$$

PROOF. Since $\nu \ll m$, and (4.8) holds $m$-a.e., it holds $\nu$-a.e. The conclusion of the corollary follows immediately, as did (4.10).

The next step in the proof is to establish (4.15) below for a $Q_\mu$ which is a convex combination of ergodic ones, that is,

$$Q_\mu = \sum_{p=1}^k \lambda_p Q\mu_p,$$

where $\lambda_p > 0, 1 \leq p \leq k$, and $\sum_{p=1}^k \lambda_p = 1$, each $Q\mu_p$ is stationary ergodic and each $\mu_p \ll m$. We will prove this case for $k = 2$. The general case can then be proven inductively. So, let

$$Q_\mu = \lambda_1 Q\mu_1 + \lambda_2 Q\mu_2, \qquad \lambda_1 > 0, \qquad \lambda_2 > 0, \qquad \lambda_1 + \lambda_2 = 1,$$

where $Q\mu_1$ and $Q\mu_2$ are stationary ergodic and $\mu_p \ll m$, $p = 1, 2$. Then

$$P^x(L_t(\omega, \cdot) \in N_{\mu,\eta})$$

$$= P^x \Big( \frac{1}{(\lambda_1 + \lambda_2)t} \int_0^{(\lambda_1 + \lambda_2)t} 1_{\cdot}(X_s(\omega)) \, ds \in N_{\mu,\eta} \Big)$$

$$(4.11) \qquad = P^x \Big( \lambda_1 \cdot \frac{1}{\lambda_1 t} \int_0^{\lambda_1 t} 1_{\cdot}(X_s(\omega)) \, ds$$

$$+ \lambda_2 \frac{1}{\lambda_2 t} \int_0^{\lambda_2 t} 1.(X_s(\theta_{\lambda_1 t} \omega)) \, ds \in N_{\mu,\eta} \Big)$$

$$\geq P^x(L_{\lambda_1 t}(\omega, \cdot) \in N_{\mu_1,\eta}, L_{\lambda_2 t}(\theta_{\lambda_1 t} \omega, \cdot) \in N_{\mu_2,\eta}),$$

by the simple convexity property of the neighborhoods. By Corollary 4.11, given $\delta > 0$, we can find a set $B$ with $\mu_1(B) > 0$ and a $t_0 > 0$ such that for all $y \in B$ and all $t \geq t_0$, we have

(4.12)     $$P^y(L_{\lambda_2 t} \in N_{\mu_2,\eta/2}) \geq \exp(-(1+\varepsilon)\lambda_2 t(H(Q\mu_2) + \delta))$$

and for $m$-a.e.$(x)$, we have, by Proposition 4.7,

(4.13)
$$\liminf_{t \to \infty} \frac{1}{t} \log E^x \Big( 1(L_{\lambda_1 t} \in N_{\mu_1,\eta}) \int_{\lambda_1 t}^{\lambda_1 t(1+\varepsilon)} 1_B(X_s) \, ds \Big)$$
$$\geq -\lambda_1(1+\varepsilon)H(Q\mu_1).$$

By (4.11), we have

(4.14)
$$P^x(L_t \in N_{\mu,\eta})$$
$$\geq E^x(1(L_{\lambda_1 t}(\omega, \cdot) \in N_{\mu_1,\eta})$$
$$\times \frac{1}{\varepsilon \lambda_1 t} \int_{\lambda_1 t}^{\lambda_1 t(1+\varepsilon)} 1_B(X_s(\omega)) \, ds \, 1(L_{\lambda_2 t}(\theta_{\lambda_1 t} \omega, \cdot) \in N_{\mu_2,\eta})),$$

hence, for all sufficiently large $t$ [cf. Remark 4.5(i)] the right-hand side is

$$\geq E^x \Big( 1(L_{\lambda_1 t}(\omega, \cdot) \in N_{\mu_1,\eta})$$
$$\times \frac{1}{\varepsilon \lambda_1 t} \int_{\lambda_1 t}^{\lambda_1 t(1+\varepsilon)} 1_B(X_s(\omega)) 1(L_{\lambda_2 t}(\theta_s \omega) \in N_{\mu_2,\eta/2}) \, ds \Big),$$

provided $\varepsilon \leq \eta \lambda_2 / 2\lambda_1$. For the last expression, we condition on $\mathcal{F}_{\lambda_1 t}^o$, take the conditioning under the integral sign and then condition on $\mathcal{F}_s^o$ first to get that it

$$= E^x \Big( 1(L_{\lambda_1 t}(\omega, \cdot) \in N_{\mu_1,\eta})$$
$$\times \frac{1}{\varepsilon \lambda_1 t} \int_{\lambda_1 t}^{\lambda_1 t(1+\varepsilon)} E^x \{ 1_B(X_s(\omega))$$
$$\times E^{X_s(\omega)} (1(L_{\lambda_2 t} \in N_{\mu_2,\eta/2})) \mid \mathcal{F}_{\lambda_1 t}^o \} \, ds \Big)$$

and now applying (4.12),

$$\geq E^x \Big( 1(L_{\lambda_1 t} \in N_{\mu_1,\eta}) \frac{1}{\varepsilon \lambda_1 t}$$

$$\times \int_{\lambda_1 t}^{\lambda_1 t(1+\varepsilon)} E^x(1_B(X_s) e^{-(1+\varepsilon)\lambda_2 t(H(Q_{\mu_2})+\delta)} \mid \mathcal{F}_{\lambda_1 t}^o) \, ds\Big)$$

$$= e^{-(1+\varepsilon)\lambda_2 t(H(Q_{\mu_2})+\delta)} E^x\Big(1(L_{\lambda_1 t} \in N_{\mu_1,\eta}) \frac{1}{\varepsilon \lambda_1 t} \int_{\lambda_1 t}^{\lambda_1 t(1+\varepsilon)} 1_B(X_s) \, ds\Big).$$

Now applying (4.13), we get via (4.14) for $m$-a.e.$(x)$ that

$$\liminf_{t\to\infty} \frac{1}{t} \log P^x(L_t \in N_{\mu,\eta}) \geq -(1+\varepsilon)(H(Q_\mu) + \delta).$$

Taking $\varepsilon_n \downarrow 0$ and $\delta_n \downarrow 0$, we can combine a countable number of $m$-null sets to conclude that for $m$-a.e.$(x)$, we have

(4.15) $$\liminf_{t\to\infty} \frac{1}{t} \log P^x(L_t \in N_{\mu,\eta}) \geq -H(Q_\mu).$$

We note that the $m$-null set may depend on the choice of $Q_\mu$ of the form $\sum_{p=1}^k \lambda_p Q\mu_p$ for given $\mu$.

We now complete the proof of Theorem 1.15. Let $Q_\mu \in \mathcal{M}_S(\Omega)$, with marginal $\mu$. Let

$$N_{\mu,\varepsilon} = \left\{ \nu \in \mathcal{M}(E) : \left| \int f_j \, d\nu - \int f_j \, d\mu \right| < \varepsilon, 1 \leq j \leq r \right\}$$

be a $\tau$-neighborhood of $\mu$. Let $g_j(\omega) = f_j(X_0(\omega))$, $1 \leq j \leq r$, and define

$$G = \left\{ Q_\nu \in \mathcal{M}_S(\Omega) : \left| \int g_j \, dQ_\nu - \int g_j \, dQ_\mu \right| < \varepsilon, 1 \leq j \leq r \right\}.$$

By Lemma 4.4, given $\delta > 0$, there exists a $Q_{\nu_0} \in G$, where $Q_{\nu_0} = \sum_{p=1}^k \lambda_p Q_{\nu_p}$, $0 < \lambda_p < 1$, $\sum \lambda_p = 1$, each $Q_{\nu_p}$ is ergodic and $|H(Q_{\nu_0}) - H(Q_\mu)| < \delta$. It follows that $\nu_0 \in N_{\mu,\varepsilon}$ and by (4.15) for $m$-a.e.$(x)$, we have

$$\liminf_{t\to\infty} \frac{1}{t} \log P^x(L_t \in N_{\mu,\varepsilon}) \geq -H(Q_{\nu_0}) \geq -(H(Q_\mu) + \delta).$$

The $m$-null set may depend on $\delta > 0$, but we can take $\delta_n \downarrow 0$ and combine the countable number of $m$-null sets to get (4.16) for $Q_\mu \in \mathcal{M}_S(\Omega)$.

We now go from an $m$-null exceptional set to an $m$-thin set. Given $Q_\mu$ with marginal $\mu$, let $N_0$ denote the $m$-null set such that if $x \notin N_0$, then (4.15) holds. If we let

$$N = \{x : R1_{N_0}(x) > 0\},$$

then $N$ is $m$-thin by definition and $x_0 \notin N$ implies that there exists $t_0$ such that $p_{t_0}(x_0, N_0^c) = 1$. Given $\delta > 0$, we can find $B \subset N_0^c$ such that $p_{t_0}(x_0, B) = \beta > 0$ and $t_1 > 0$ such that for all $t \geq t_1$,

(4.16) $$\inf_{y \in B} \frac{1}{t} \log P^y(L_t \in N_{\mu,\varepsilon}) \geq -H(Q_\mu) - \delta.$$

Now, using Remark 4.5(i), we have for all $t$ sufficiently large (depending on $t_0$ and $\varepsilon$),

$$\frac{1}{\beta} P^{x_0}(L_t \in N_{\mu,2\varepsilon}) \geq \frac{1}{\beta} \int_B p_{t_0}(x_0, dy) P^y(L_t \in N_{\mu,\varepsilon}).$$

By Jensen's inequality,

$$-\log \beta + \log P^{x_0}(L_t \in N_{\mu,2\varepsilon}) \geq \frac{1}{\beta} \int_B p_{t_0}(x_0, dy) \log P^y(L_t \in N_{\mu,\varepsilon})$$

and now dividing by $t$ and letting $t \to \infty$, because of the lower boundedness [by (4.16)], we apply Fatou's lemma to get

$$\liminf_{t \to \infty} \frac{1}{t} \log P^{x_0}(L_t \in N_{\mu,2\varepsilon}) \geq -H(Q_\mu) - \delta.$$

$\delta > 0$ being arbitrary, it can be dropped on the right-hand side and we get (4.15) for $x_0 \notin N$. For a given $\mu \in \mathcal{M}(E)$, $Q_\mu$ is by no means unique and $N$ may depend on the choice of $Q_\mu$. However, by the contraction principle, there exists a sequence $Q_\mu^{(n)}$ such that $H(Q_\mu^{(n)}) \to I(\mu)$ as $n \to \infty$. The countable number of $m$-thin sets corresponding to each $Q_\mu^{(n)}$ may be combined and we finally conclude that if $\mu \ll m$, then there exists an $m$-thin set $N$ such that if $x \notin N$, then

$$\liminf_{t \to \infty} \frac{1}{t} \log P^x(L_t \in N_{\mu,\varepsilon}) \geq -I(\mu).$$

If Condition 1.5 holds with the associated framework, then the $m$-thin set has capacity 0 and Theorem 1.15 is established.  $\square$

REMARK 4.17.   If the semigroup is Feller, since a weak neighborhood is a $\tau$-neighborhood, we get a stronger result than the one which applies only to $w$-neighborhoods.

PROOF OF COROLLARY 1.16.   Let $\alpha = \inf\{I(\mu) : \mu \in U\}$. We may assume that $\alpha < \infty$. For $j \geq 1$, there exist $\mu_j \in U$ such that $I(\mu_j) < \alpha + j^{-1}$. Let $V_j$ be a $\tau$-neighborhood of $\mu_j$ contained in $U$. Then, by Theorem 1.15, there exists an $m$-thin set $N_j$ such that if $x \notin N_j$, then

$$\liminf_{t \to \infty} \frac{1}{t} \log P^x(L_t \in V_j) \geq -(\alpha + j^{-1}).$$

It follows that if $x \notin N = \bigcup_{j=1}^{\infty} N_j$, then

$$\liminf_{t \to \infty} \frac{1}{t} \log P^x(L_t \in U) \geq -(\alpha + j^{-1})$$

for all $j \geq 1$, hence the corollary follows.  $\square$

**5. The self-adjoint case.** In this section, we assume that $\{T_t\}$ is a symmetric semigroup, so its $L_2(m)$-generator $L$ is self-adjoint, the domain $\mathcal{D}(L)$ of which is dense in $L_2(m)$. We denote by $H$ the canonical square root of $-L$. Then $\mathcal{D}(L) \subset \mathcal{D}(H)$. As far as the Markov process is concerned, we only need it to be progressively measurable.

To prove Theorem 1.17, we need some lemmas. Let $v \in D_1$ and $\varepsilon > 0$ be fixed. With $v_\varepsilon := v + \varepsilon$, let

$$V = \hat{L}v / v_\varepsilon.$$

Recall that $V$ is a bounded $L_1(m)$ function. If $f$ is a function which is in an $L_2(m)$-equivalence class, define

$$(5.1) \qquad \tilde{p}_t f(x) = E^x(f(X_t)e^{-\int_0^t V(X_s)ds}).$$

Let $\tilde{T}_t f$ denote the corresponding operator on $L_2(m)$.

LEMMA 5.2.  $\{\tilde{T}_t, t \geq 0\}$ *is a self-adjoint strongly continuous semigroup on $L_2(m)$, its $L_2(m)$-generator $\tilde{L}$ is given by*

$$(5.3) \qquad \tilde{L}f = Lf - Vf$$

*and $\mathcal{D}(\tilde{L}) = \mathcal{D}(L)$.*

PROOF.  The fact that $\{\tilde{T}_t\}$ is a strongly continuous self-adjoint semigroup in $L_2(m)$ is proved in [2], page 130. The function $V$ is bounded and in $L_2(m)$, hence, for any function $f$ in $L_2(m)$,

$$(p_t f - \tilde{p}_t f)(x) = E^x(f(X_t) - f(X_t)e^{-\int_0^t V(X_s)ds})$$
$$= E^x\left(f(X_t)\int_0^t V(X_s)\,ds + f(X_t)O(t^2)\right)$$

as $t \to 0$, where $O(t^2)$ is uniform in $\omega$. Therefore,

$$(5.4) \qquad \begin{aligned} &\limsup_{t\to 0}\left\|\frac{1}{t}(p_t f - \tilde{p}_t f) - Vf\right\|_2 \\ &= \limsup_{t\to 0}\left\|E^x\left(f(X_t)\cdot\frac{1}{t}\int_0^t V(X_s)\,ds\right) - V(x)f(x)\right\|_2.\end{aligned}$$

Letting $W(t) = t^{-1}\int_0^t V(X_s)\,ds$, the right-hand side in (5.4) is

$$\leq \limsup_{t\to 0}\|E^x((f(X_t) - f(x))W(t))\|_2$$
$$+ \limsup_{t\to 0}\|f(x)E^x W(t) - f(x)V(x)\|_2.$$

Since $|W(t)| \le c$ for some $c > 0$ independent of $t$,

$$\|E^x f(X_t) W(t) - f(x) E^x W(t)\|_2^2 \le c^2 \int (E^x f^2(X_t) - 2E^x f(X_t) f + f^2) \, dm$$

$$= c^2 \left( 2\|f\|_2^2 - 2 \int (p_t f) f \, dm \right)$$

$$\to 0 \qquad \text{as } t \to 0.$$

Also,

$$\|f(x) E^x W(t) - f(x) V(x)\|_2^2 = \int f^2(x) (E^x (W(t) - V(x)))^2 \, dm(x).$$

$W$ and $V$ are bounded and $\int (E^x (W(t) - V(x)))^2 \, dm(x) \to 0$ as $t \to 0$ so that $E^x (W(t) - V(x)) \to 0$ in $m$-measure and by dominated convergence, $\|f E^x W(t) - f V\|_2 \to 0$ as $t \to 0$. Therefore, the right-hand side in (5.4) is 0. It follows that for any $f \in L_2(m)$,

$$\left\| \frac{T_t f - f}{t} - \frac{\tilde{T}_t f - f}{t} - V f \right\|_2 \to 0 \qquad \text{as } t \to 0,$$

hence, $f \in \mathcal{D}(L)$ if and only if $f \in \mathcal{D}(\tilde{L})$ and $\tilde{L} f = L f - V f$, which proves the lemma.  $\square$

LEMMA 5.5.   $\{\tilde{T}_t\}$ is a contraction on $L_2(m)$.

PROOF.   Let $v, g \in D_0$. Let $u_n \in D$, $u_n \nearrow 1$ and for some $\eta > 0$, let

$$v_n = \frac{1}{\eta} \int_0^\eta p_s u_n \, ds.$$

Then $v + \varepsilon v_n \nearrow v_\varepsilon (= v + \varepsilon)$ and by the self-adjointness of $\{\tilde{T}_t\}$ on $L_2(m)$ and the Feynman–Kac formula (Corollary 2.6), we have that for $n \ge 1$,

$$\int (v + \varepsilon v_n) \tilde{p}_t g \, dm = \int \tilde{p}_t (v + \varepsilon v_n) g \, dm = \int (v + \varepsilon v_n) g \, dm.$$

Letting $n \to \infty$, since $0 \le v + \varepsilon v_n \le v_\varepsilon$ and $\tilde{p}_t g$ and $g$ belong to $L_1(m)$, we get

$$(5.6) \qquad \int v_\varepsilon \tilde{p}_t g \, dm = \int v_\varepsilon g \, dm.$$

By Lemma 2.3, $D_0$ is dense in $L_1^+(m)$, hence (5.6) holds for all $g \in L_1^+(m)$. This shows that $v_\varepsilon \, dm$ is an invariant measure for $\tilde{p}_t$.

If $f \in D_0$, then

$$\int |\tilde{T}_t f|^2 \, dm = \int \{E^x (f(X_t) e^{-\int_0^t V(X_s) \, ds})\}^2 \, dm$$

and by the Schwarz inequality,

$$\leq \int E^x \{v_\varepsilon(X_t))^{1/2} e^{-(1/2) \int_0^t V(X_s) ds}\}^2$$

$$\times E^x \{f(X_t)(v_\varepsilon(X_t))^{-1/2} e^{-(1/2) \int_0^t V(x_s) ds}\}^2 \, dm$$

$$= \int E^x (v_\varepsilon(X_t) e^{-\int_0^t V(X_s) ds}) E^x (f^2(X_t)(v_\varepsilon(X_t))^{-1} e^{-\int_0^t V(X_s) ds}) \, dm.$$

By the Feynman–Kac formula (Corollary 2.6),

$$= \int v_\varepsilon(x) \tilde{T}_t(f^2/v_\varepsilon)(x) \, dm(x)$$

and by (5.6) applied to $g = f^2/v_\varepsilon$ (which is in $L^1$ if $f \in D_0$), we get that the last expression

$$= \int v_\varepsilon(f^2/v_\varepsilon) \, dm = \int f^2 \, dm.$$

Since $D_0$ is dense in $L_2^+(m)$ and $\{\tilde{T}_t\}$ is strongly continuous on $L_2(m)$, we have shown that

$$\|\tilde{T}_t f\|_2 \leq \|f\|_2$$

for all $f \in L_2^+(m)$, which suffices, and the lemma is proved. $\square$

LEMMA 5.7. *Let $\mu \ll m$ and $f^2 = d\mu/dm$, where $f$ as a function represents a nonnegative version of the square root of $f^2$. Then, for all $h \geq 0$,*

$$\int (f - p_h f) f \, dm \leq h I(\mu).$$

PROOF. Let $f_n = f \wedge n$ for $n \geq 1$, and let $t_k \downarrow 0$ and

$$f_n^{(k)} = \frac{1}{t_k} \int_0^{t_k} p_s f_n \, ds$$

for $k \geq 1$. Then $f_n^{(k)} \in D_0$ and $f_n^{(k)} \leq n$. By (2.18), if $v \in D_0, \varepsilon > 0$ and $v_\varepsilon = v + \varepsilon$, we have

$$\int \log \frac{p_h v_\varepsilon}{v_\varepsilon} \, d\mu \geq -h I(\mu).$$

Since $\log x \leq x - 1$ for $x > 0$, this gives

$$\int \frac{v_\varepsilon - p_h v_\varepsilon}{v_\varepsilon} \, d\mu \leq h I(\mu).$$

Taking $v = f_n^{(k)}$, we then get

$$\int \frac{f_n^{(k)} - p_h f_n^k}{f_n^{(k)} + \varepsilon} f^2 \, dm \leq h I(\mu).$$

As $k \to \infty$, we get, by dominated convergence,

$$\int \frac{f_n - p_h f_n}{f_n + \varepsilon} f^2 \, dm \le h I(\mu).$$

Since

$$\frac{p_h f_n}{f_n + \varepsilon} \le \frac{n}{n + \varepsilon} 1_{[f > n]} + \frac{p_h f}{f + \varepsilon} 1_{[f \le n]},$$

we can let $n \to \infty$ to get

$$\int \frac{f - p_h f}{f + \varepsilon} f^2 \, dm \le h I(\mu).$$

Finally, $|f - p_h f| f^2/(f + \varepsilon) \le |f - p_h f| f \in L_1(m)$, so we can let $\varepsilon \downarrow 0$ and the lemma is proved. □

PROOF OF THEOREM 1.17. Let $\mu \ll m, d\mu = f^2 \, dm$ and $f \ge 0$. By Lemma 5.7, we have that for $h > 0$,

$$\frac{1}{h} \int (f - T_h f) f \, dm \le I(\mu).$$

If $\{E_\lambda, \lambda \ge 0\}$ is the spectral measure for the positive definite operator $-L$, then following the argument of Donsker and Varadhan ([3], page 46), this implies that

$$\int_0^\infty \lambda d \langle E_\lambda f, f \rangle \le I(\mu),$$

so $I(\mu) < \infty$ implies that $f \in \mathcal{D}(H)$ and

$$\|Hf\|_2^2 \le I(\mu).$$

For the other direction, let $v \in D_1, \varepsilon > 0$ and $v_{\bar{\varepsilon}} = v + \varepsilon$. Let $V = \hat{L} v / v_\varepsilon$ and $\tilde{T}_t$ be given via (5.1). Since, by Lemma 5.5, $\tilde{T}_t$ is a contraction on $L_2(m)$, we have

$$\langle \tilde{T}_t f - f, f \rangle \le 0,$$

therefore, if $f \in \mathcal{D}(L)$, we get

$$\langle \tilde{L} f, f \rangle \le 0,$$

which gives [by (5.3)]

$$(5.8) \qquad \langle -Lf + Vf, f \rangle \ge 0,$$

meaning that for every $v \in D_1$, $f \in \mathcal{D}(L)$ and $\varepsilon > 0$,

$$\langle Lf, f \rangle \le \langle Vf, f \rangle,$$

that is,

$$\int \frac{Lv}{v_\varepsilon} f^2 \, dm \geq \langle Lf, f \rangle.$$

Therefore,

$$\inf_{\substack{v \in D_1 \\ \varepsilon > 0}} \int \frac{Lv}{v_\varepsilon} f^2 \, dm \geq \langle Lf, f \rangle.$$

It follows that

$$(5.9) \qquad I(\mu) \leq -\langle Lf, f \rangle = \|Hf\|_2^2$$

if $f \in \mathcal{D}(L)$. If $f \in \mathcal{D}(H)$, then we take $0 < \lambda_n \nearrow \infty$ and $f_n := E_{\lambda_n} f$. From the spectral representation, it is clear that $f_n \in \mathcal{D}(L)$, $f_n \to f$ in $L_2(m)$ and $Hf_n \to Hf$. Since (5.9) holds for each $f_n$, letting $n \to \infty$, we then see that for $f \in \mathcal{D}(H)$,

$$\liminf_{n \to \infty} I(\mu_n) \leq \|Hf\|_2^2,$$

where $d\mu_n = f_n^2 \, dm$. Since $I$ is lower semicontinuous with respect to the $\tau$-topology on $\mathcal{M}(E)$, we conclude that $I(\mu) \leq \liminf_{n \to \infty} I(\mu_n)$, where $d\mu = f^2 \, dm$, and the theorem is proved. $\square$

REMARK 5.10. It is easy to check that if $\mathcal{E}(u, v)$ is the associated Dirichlet form, then $u \in \mathcal{D}(H)$ if and only if $u \in \mathcal{D}(\mathcal{E})$, and if $u \in \mathcal{D}(H)$, then

$$\|Hu\|_2^2 = \mathcal{E}(u, u)$$

(cf. [7], page 28).

REMARK 5.11 (Relation between Theorem 1.17 and Theorem 5 of [3]). Donsker and Varadhan [3] consider a semigroup $\{T_t, t \geq 0\}$ acting on $B_b(E)$. Let

$$B_0 = \left\{ u \in B_b(E) : \lim_{t \to 0} \|p_t u - u\|_\infty = 0 \right\},$$

where $\|\cdot\|_\infty$ denotes the sup norm. They assume that $m$ is a $\sigma$-finite invariant measure for the semigroup and that $B_0 \cap L_2(m)$ is dense in $L_2(m)$. Then, for $u \in L_2(m)$, given $\varepsilon > 0$, there exists $u_0 \in B_0 \cap L_2(m)$ such that $\int |u - u_0|^2 \, dm < \varepsilon$. Furthermore,

$$\int |p_t u - u|^2 \, dm \leq 3 \left\{ \int |p_t u - p_t u_0|^2 \, dm \right.$$
$$\left. + \int |p_t u_0 - u_0|^2 \, dm + \int |u_0 - u|^2 \, dm \right\}$$
$$\leq 3 \left\{ 2 \int |u - u_0|^2 \, dm + \int |p_t u_0 - u_0|^2 \, dm \right\}.$$

Therefore, to show that $\{T_t\}$ is a strongly continuous semigroup on $L_2(m)$, it suffices to show that for $u_0 \in B_0 \cap L_2(m)$,

$$\lim_{t \to 0} \int |p_t u_0 - u_0|^2 \, dm = 0.$$

We have

$$\int |p_t u_0 - u_0|^2 \, dm = \int (p_t u_0)^2 \, dm - 2 \int (p_{t/2} u_0)^2 \, dm + \int u_0^2 \, dm$$

$$\leq \int p_t u_0^2 \, dm - 2 \int (p_{t/2} u_0)^2 \, dm + \int u_0^2 \, dm$$

and since $m$ is invariant for the semigroup, the last expression

$$= 2 \int u_0^2 \, dm - 2 \int (p_{t/2} u_0)^2 \, dm.$$

Since $u_0 \in B_0$, by Fatou's lemma, we have

$$\liminf_{t \to 0} \int (p_{t/2} u_0)^2 \, dm \geq \int u_0^2 \, dm,$$

hence, $\limsup_{t \to 0} \int |p_t u_0 - u_0|^2 \, dm = 0$, the semigroup $\{T_t\}$ is strongly continuous on $L_2(m)$ and our results apply without some extraneous assumptions.

**6. The discrete-time case.** We assume that $\{T_n, n \geq 1\}$ is a semigroup on $L_2(m)$ for which (1.2) holds. Let $\hat{I}_1(\mu)$ be as defined in Section 1. Let

$$I_1(\mu) = \begin{cases} \hat{I}_1(\mu), & \text{if } \mu \ll m, \\ \infty, & \text{otherwise.} \end{cases}$$

Then, if one uses the discrete analog of the Feynman–Kac formula (cf. [3] or [2]), a proof analogous to that of Theorem 1.8 yields the following:

THEOREM 6.1. *Let $\mu \in \mathcal{M}(E)$ and let $a < \hat{I}_1(\mu)$. There then exists a $\tau$-neighborhood $N_\mu$ of $\mu$ such that*

$$\limsup_{n \to \infty} \frac{1}{n} \log \sup_{x \in E} P^x(L_n \in N_\mu) \leq -a.$$

REMARK 6.2. This immediately implies the upper bound for $\tau$-compact sets, as in the continuous-time case. Furthermore, if the semigroup is Feller, then $N_\mu$ may be chosen to be a weak neighborhood in the statement of the theorem and "$\tau$-compact" can be replaced by "$w$-compact" in the corollary.

The next theorem is the analog of Theorem 1.11, but first, we define an $m$-thin set in the present context.

A set $N$ is $m$-thin in the discrete-time context if there exists a set $B$ such that $m(B) = 0$ and

$$N = \left\{ x \in E : \sum_{n=1}^{\infty} 2^{-n} p_n(x, B) > 0 \right\}.$$

We then have the following.

THEOREM 6.3. *Let $\mu \in \mathcal{M}(E)$ and let $a < I_1(\mu)$. There then exists a $\tau$-neighborhood $N_\mu$ of $\mu$ and an $m$-thin set $N$ such that if $x \notin N$, then*

$$\limsup_{n \to \infty} \frac{1}{n} \log P^x(L_n \in N_\mu) \leq -a.$$

For the lower bound, we also need the counterpart of Condition 1.14, which reads as follows.

(6.4)  If $A \in \mathcal{B}(E)$ and $m(A) > 0$,    then for $m$-a.e.$(x)$, $\displaystyle\sum_{n=1}^{\infty} p_n(x, A) > 0$.

Under this additional assumption, we have the following.

THEOREM 6.5. *Let $\mu \in \mathcal{M}(E)$ and let $N_\mu$ be a $\tau$-neighborhood of $\mu$. Then there exists an $m$-thin set $N$ such that if $x \notin N$,*

$$\liminf_{n \to \infty} \frac{1}{n} \log P^x(L_n \in N_\mu) \geq -I(\mu).$$

REMARK 6.6. The analog of Corollary 1.16 also holds for $U$ a $w$-open subset of $\mathcal{M}(E)$.

**Acknowledgments.** We are indebted to the referees for their careful reading of the manuscript and for making several useful suggestions. We have followed their suggestions as much as possible and hope that the revised version is easier to understand and free of misprints. We are also grateful to one of them for drawing attention to the work of Takeda [10]. His definition of the rate function is similar to ours and is given directly in terms of the resolvent operators associated with the $L_2$-semigroup; our definition is given in terms of the generator. We should note that all of these definitions are essentially similar to the ones given by Donsker and Varadhan in their papers.

We thank Michael Röckner for bringing the works in [7] and [9] to our attention.

DEPARTMENT OF MATHEMATICS
UNIVERSITY OF MINNESOTA
MINNEAPOLIS, MINNESOTA 55455
USA
E-MAIL: jain@math.umn.edu
          krylov@math.umn.edu